\documentclass{amsart}
\usepackage{amsfonts}
\usepackage{enumerate}
\usepackage{amsthm}
\usepackage{latexsym}
\usepackage{amsmath}

\theoremstyle{definition}
\newtheorem{defi}{Definition}[section]
\newtheorem{remark}{Remark}
\newtheorem{example}{Example}
\theoremstyle{plane}
\newtheorem{theo}{Theorem} 
\newtheorem{lem}{Lemma}[section]
{}
\newtheorem{proposition}{Proposition}[section]
\newtheorem{cor}{Corollary}[section]

\newcommand{\R}{{\mathbb R}}
\newcommand{\Z}{{\mathbb Z}}
\newcommand{\N}{{\mathbb N}}
\newcommand{\con}{{\rm conv}}
\newcommand{\al}{\mathbf{1}}
\newcommand{\Ehr}{{\rm Ehr}}

\title{Gorenstein polytopes obtained from bipartite graphs}
\author{Makoto Tagami}
\address{Universit\"at Magdeburg, Institut f\"ur Algebra und Geometrie,\\
  Universit\"ats\-platz 2, D-39106 Magdeburg, Germany}
\email{tagami@kenroku.kanazawa-u.ac.jp}
\thanks{The author was supported by Research Fellowships of the Japan Society for the Promotion of Science for Young Scientists.}
\begin{document}
\begin{abstract}
Beck et. al. characterized the grid graphs whose perfect matching polytopes are Gorenstein and they also showed that for some parameters, perfect matching polytopes of torus graphs are Gorenstein. In this paper, we complement their result, that is, we characterize the torus graphs whose perfect matching polytopes are Gorenstein. Beck et. al. also gave a method to construct an infinite family of Gorenstein polytopes.  In this paper, we introduce a new class of polytopes obtained from graphs and we extend their method to construct many more Gorenstein polytopes.
\end{abstract}
\maketitle
\section{Introduction}\label{intro}

Lattice polytopes are polytopes whose vertices all are lattice points. For $S \subset \R^n$ and $t\in \N$, we put $tS=\{tx \mid x\in S\}$ and $L_S(t)=\sharp(tS\cap \Z^n)$. Ehrhart \cite{ehrhart} proved that for a $d$-dimensional lattice polytope $P$, $L_P(t)$ is always a polynomial of degree $d$ in $t$. $L_P(t)$ is called the \textit{Ehrhart polynomial} of $P$. Also the formal power series $\Ehr_P(z)=1+\sum_{t\in \N} L_P(t)z^t$ is called the \textit{Ehrhart series} of $P$. Since $L_P(t)$ is a polynomial of degree $d$, the Ehrhart series of $P$ can be written as the rational function: 
\[\Ehr_P(z)=\frac{\sum_{i=0}^s h_iz^i}{(1-z)^{d+1}},\]
where $s\le d$. $s$ and $r=d+1-s$ are  called the \textit{degree} and \textit{codegree} of $P$, respectively. The polynomial of the numerator is called the $h^*$-\textit{polynomial} of $P$. It is well-known that $h_0=1$, and the codegree $r$ is equal to the minimal integer $t$ for which $tP^\circ$ contains a lattice point, and $h_s=\sharp(rP^\circ\cap Z^n)$. Here, for $S\subset \R^n$, $S^\circ$ denotes the relative interior of $S$. As a general reference on the Ehrhart theory of lattice polytopes we refer to the recent book of Matthias Beck and Sinai Robins \cite{Beck} and the references within.

Ehrhart polynomials have also an algebraic meaning in the sence that the Ehrhart polynomial of a polytope $P$ can be interpreted as the Hilbert function of the Ehrhart ring of $P$. We say $P$ is \textit{Gorenstein} when the Ehrhart ring is Gorenstein(see \cite{Hibi} and \cite{Stanley} for Ehrhart rings and Gorenstein property). $P$ is Gorenstein if and only if the coefficients of the $h^*$-polynomial are symmetric, that is, $h_i=h_{s-i}$ for any $i$. In terms of Ehrhart polynomials, this is equivalent to $L_{P^\circ}(r)=1$ and $L_P(t-r)=L_{P^\circ}(t)$ for any $t>r$. 

Let $G=(V,E)$ be an undirected graph without multiple edges and loops. In fact, even if there are multiple edges, the argument below also holds after a little change. Here $V$ and $E$ denote the vertex set and the edge set, respectively. $M\subset E$ is a \textit{matching} if any two distinct edges do not intersect. 
If every vertex lies on some edge in $M$, we call $M$ a \textit{perfect matching} for $G$. For a perfect matching $M$, we define the characteristic vector $\chi_M \in \R^E$ as follows: for $e\in E$,
\[(\chi_M)_e:=
\begin{cases} 
1: & \mbox{if $e \in M$},\\
0: & \mbox{otherwise}.
  \end{cases}\]
\begin{defi}[Perfect matching polytope]
The \textit{perfect matching polytope} of $G$, $P_G$ is defined to be the convex hull in $\R^E$ of the characteristic vectors of all perfect matchings:
\[P_G:=\con\{\chi_M\mid M :\,\mbox{a perfect matching of $G$} \} \subset \R^E.\]
\end{defi}
In general, $P_G$ is not full-dimensional. We note that interior lattice points of $P_G$ mean lattice points in the relative interior of $P_G$. 
By Edmond's famous theorem we know a hyperplane description for perfect matching polytopes:
\begin{theo}[Edmond \cite{edmond}]
Let $G=(V,E)$ be a graph with $|V|$ even. Then $x=(x_e)_{e\in E}\in \R^E$ lies in $P_G$ if and only if the following conditions hold:
\begin{enumerate}[$(1)$]
\item $x_e\ge 0 \;(\forall e \in E)$,
\item $\sum_{v \in e}x_e=1$ $(\forall v\in V)$,
\item $\sum_{e\in C(S,S')}x_e \ge 1$ ($\forall S \subset V$, $|S|$ is odd),
\end{enumerate}
where $v \in e$ means that $v$ is incident to $e$, and $S'$ denotes the complement set of $S$ in $V$, and for subsets $S$ and $T$ of $V$, $C(S,T):=\{(u,v)\in E\mid u\in S,\;v\in T\}$. 
\end{theo}
A graph $G=(V,E)$ is \textit{bipartite} if there exists some partition $V=V_1\cup V_2$ such that $C(V_i,V_i)=\emptyset$ for $i=1,2$.
It is well-known that, if a graph is bipartite, then we can omit the third condition, that is, 
$x\in P_G$ if and only if the conditions (1) and (2) hold. For a subset $S$ of $V$, we call edges in $C(S,S')$ \textit{bridges} from $S$. We refer to Gr\"{o}tschel-Lov\'{a}sz-Schrijver \cite{GLS} about perfect matching polytopes.  

The $m\times n$ grid graph $\mathcal{G}(m,n)=(V,E)$ is defined as follows: $V:=\{(i,j)\mid 0\le i \le m-1,\,0\le j \le n-1\}$, and $((i,j),(k,l))\in E$ if and only if $|i-k|+|j-l|=1$.
The $m\times n$ torus graph $\mathcal{G}_T(m,n)$ consists of the same vertex and edge set as $\mathcal{G}(m,n)$ with 
the additional edges $\{((0,j),(m-1,j))\mid 0\le j \le n-1\}$ and $\{((i,0),(i,n-1))\mid 0\le i\le m-1\}$.
 
Using Edmond's theorem, Beck et. al \cite{BHS} characterized the grid graphs whose perfect matching polytopes are Gorenstein. They also showed  that the perfect matching polytopes of torus graphs for some parameters are Gorenstein. We denote by $\mathcal{P}(m,n)$ and $\mathcal{P}_T(m,n)$ the perfect matching polytopes of $\mathcal{G}(m,n)$ and $\mathcal{G}_T(m,n)$, respectively. That is to say that, they showed the following:
\begin{theo}[B-H-S \cite{BHS}]\label{bhs}
If $m=1$ or $m$ is even, and $n$ is even, then $\mathcal{P}_T(m,n)$ is Gorenstein.
\end{theo}
In section \ref{torus}, we complement Theorem \ref{bhs} by showing:
\begin{theo}\label{hokan}
If $mn$ is even, then $\mathcal{P}_T(m,n)$ is Gorenstein if and only if $m=1$ or even, and $n$ is even, or $(m,n)=(2,3),\, (2,5)$. 
\end{theo}
We remark that Beck et. al. \cite{BHS} claimed Theorem \ref{bhs} is a corollary of the more general result:
\begin{theo}[B-H-S \cite{BHS}]\label{bhs3}
Let $G$ be a $k$-regular bipartite graph with even vertices. Then the perfect matching polytope $P_G$ is Gorenstein.  
\end{theo}
Here a graph $G$ is $k$-\textit{regular} if any vertex is incident to exactly $k$ edges. 
Theorem \ref{bhs3} constructs an infinite family of Gorenstein polytopes. 
In section \ref{smp} we introduce a new class of polytopes obtained from graphs which are a natural extension of perfect matching polytopes. For these polytopes we show an analogous result to Edmond's theorem. Also using these polytopes, we extend Theorem \ref{bhs3} in order to construct many more Gorenstein polytopes. For an another method to construct Gorenstein polytopes from graphs, we refer to Ohsugi-Hibi \cite{hibi2}.

\section{Torus graphs and Perfect matching polytopes}\label{torus}
In this section, we complement the characterization for torus graphs whose perfect matching polytope is Gorenstein.
\begin{lem}\label{first}
Let $\mathcal{G}_T(m,n)=(V,E)$ be an $m\times n$ torus graph and $m,\,n \ge 3$. Then for any subset $S \subset V\;(2\le |S|\le |V|-2)$, there are at least 6 bridges.
\end{lem}
\begin{proof}
We call points of $S$ black points and points of $S'$ white points, respectively. If, for any column all points on the column are  black or all points  are white, then the Lemma follows since, in that case, there are at lease 2 bridges on each row and $m\ge 3$. Threrefore we may assume that there exists some column on which there are both black points and white points. 
Without loss of generality we may assume that such a column is the first column and that $(1,1)$ is white. Already there are at least 2 bridges on the first column. We divide the cases into (I) the case when there is a black point on the first row, and (II) the other case.

(I) In this case, there are at least 2 bridges on the first row. So we have $4$ bridges already. 
 Since $2\le |S| \le |V|-2$, there is a white point except for $(1,1)$. Put such a point to be $(i,j)$. Without loss of generality, we may assume  $i\not=1$. If there exists a black point on the $i$-th row, then the number of bridges increases by at least 2, and so the Lemma follows. Next we assume that all points on the $i$-th row are white.  We let a black point on the first row be $(1,k)$, then there exist at least $2$ bridges on the $k$-th column. Therefore in this case, we have at least $6$ bridges. 

(II) In this case, if there exists a black point on each column, then since there are at least 2 bridges on each column and $n\ge 3$, the Lemma follows. Therefore we may assume that there exists some column such that all points on the column are white. In this case a black point on the first column has a white point on both the row and the column on which the black point lies. So exchanging the position of black points and white points we have the same situation as in (I).
\end{proof}
Next we show a lemma about interior lattice points in perfect matching polytopes.
\begin{lem}\label{rel}
Let $P$ be a polytope defined by conditions $(1)$, $(2)$ and $(3)$ in Edmond's theorem. Set
\begin{enumerate}[$(1')$]
\item $x_e>0 \;(\forall e \in E),$ 
\item $\sum_{v \in e}x_e=1\; (\forall v\in V),$ 
\item $\sum_{e\in C(S,S')}x_e >1\; (\forall S \subset V,\, 3\le |S|\le |V|-3$\, odd). 
\end{enumerate}
Assume that there exists a vector $x$ satisfying $(1')$, $(2')$ and $(3')$. 
Then $x\in P^\circ$ and the relative interior of $P$ is given by these conditions $(1')$, $(2')$ and $(3')$.
\end{lem}
\begin{proof}
Take a vector $x$ satisfying $(1')$, (2') and (3'). Denote by $W$ the linear subspace defined by equations $\sum_{v \in e}x_e=0\; (\forall v\in V)$. Then $P$ lies on the affine subspace $x+W$. If the norm of $y\in W$ is enough small, then $x+y$ also satisfies (1'), (2') and (3'). This implies $\dim P=\dim W$ and $x\in P^\circ$ in the sence of the relative interior. After all, we see that the relative interior of $P$, $P^\circ$ is defined by (1'), (2') and (3').  
\end{proof}
Let conditions $(1)_t$, $(2)_t$ and $(3)_t$ be 
\begin{enumerate}[$(1)_t$]
\item $x_e\ge 0 \;(\forall e \in E),$ 
\item $\sum_{v \in e}x_e=t\; (\forall v\in V),$ 
\item $\sum_{e\in C(S,S')}x_e \ge t\; (\forall S \subset V,\, 3\le |S|\le |V|-3$\, odd), 
\end{enumerate}
 and conditions  $(1')_t$, $(2')_t$ and $(3')_t$ be 
\begin{enumerate}[$(1')_t$]
\item $x_e>0 \;(\forall e \in E),$ 
\item $\sum_{v \in e}x_e=t\; (\forall v\in V),$ 
\item $\sum_{e\in C(S,S')}x_e >t\; (\forall S \subset V,\, 3\le |S|\le |V|-3$\, odd). 
\end{enumerate}
Denote all-one vector $(1,\ldots,1)$ by $\al$ and define $\iota:\R^E\longrightarrow \R^E$ by $\iota(x)=x+\al$.
\begin{lem}\label{inj}
Assume that $\al$ satisfies conditions $(1')_k$, $(2')_k$ and $(3')_k$. Then $\iota$ gives an injective map from $lP\cap \Z^E$ to $(l+k)P^\circ\cap \Z^E$
\end{lem}
\begin{proof}
Since $\al$ satisfies $(1')_k$, $(2')_k$ and $(3')_k$, by Lemma \ref{rel}, $tP^\circ$ is defined by $(1')_t$, $(2')_t$ and $(3')_t$.
Let $x \in lP \cap \Z^E$. Then $x$ satisfies conditions $(1')_l$, $(2')_l$ and $(3')_l$. Since $\al$ satisfies conditions $(1')_k$, $(2')_k$ and $(3')_k$, by summing two corresponding inequalities for $x$ and $\al$, we see that $\iota(x)=x+\al$ satisfies $(1')_{l+k}$, $(2')_{l+k}$ and $(3')_{l+k}$. Clearly $\iota(x)\in \Z^E$. Therefore $\iota(x)\in (l+k)P^\circ\cap \Z^E$. 
\end{proof}
We know the dimension of perfect matching polytopes of grid graphs and torus graphs.
\begin{proposition}[B-H-S \cite{BHS}]\label{bhs2}
If $mn$ is even. Then 
\begin{enumerate}[$(i)$]
\item $\dim \mathcal{P}(m,n)=(m-1)(n-1)$,
\item if $n>2$ is even, then $\dim \mathcal{P}_T(2,n)=n+1$,
\item if $m>2$ and $n>2$ are both even, then $\dim \mathcal{P}_T(m,n)=mn+1$,
\item if $n>1$ is odd, then $\dim \mathcal{P}_T(2,n)=n$,
\item if $m>2$ is even and $n=1$, then $\dim \mathcal{P}_T(m,n)=1$,
\item if $m>2$ is even and $n>1$, then $\dim \mathcal{P}_T(m,n)=mn$.
\end{enumerate}
\end{proposition}
Now we can show Theorem \ref{hokan}.
\begin{proof}[Proof of Theorem \ref{hokan}]
Let $P=\mathcal{P}_T(m,n)$. First we show the sufficiency. 
The case when $m=1$ or $m$ is even, and $n$ is even  has been shown in Theorem \ref{bhs}. 

Let $(m,n)=(2,3)$. We show that $L_{P^\circ}(3)=0$ and $L_{P^\circ}(4)\not=0$. Define $x\in \R^E$ as in Figure 1, then $x$ satisfies conditions $(1')_4$, $(2')_4$ and $(3')_4$. So by Lemma \ref{rel}, we see that $x\in 4P^\circ$ and $tP^\circ$ is defined by $(1')_t$, $(2')_t$ and $(3')_t$. Since the graph is $3$-regular, lattice points satisfying $(1')_3$ and $(2')_3$ must be assigned by 1 on each edge. But if we take $S=\{(0,0),\,(0,1),\,(0,2)\}$, then the condition $(3')_3$ does not hold. Hence $3P^\circ$ has no lattice points. 
\begin{figure}[htbp]
\begin{center}
\unitlength 0.1in
\begin{picture}( 13.4500,  7.3500)( 14.5500,-20.8000)
%
\special{pn 20}%
\special{sh 1}%
\special{ar 1800 1600 10 10 0  6.28318530717959E+0000}%
\special{sh 1}%
\special{ar 2200 1600 10 10 0  6.28318530717959E+0000}%
\special{sh 1}%
\special{ar 2600 1600 10 10 0  6.28318530717959E+0000}%
%
\special{pn 20}%
\special{sh 1}%
\special{ar 1800 2000 10 10 0  6.28318530717959E+0000}%
\special{sh 1}%
\special{ar 2200 2000 10 10 0  6.28318530717959E+0000}%
\special{sh 1}%
\special{ar 2600 2000 10 10 0  6.28318530717959E+0000}%
%
\special{pn 8}%
\special{pa 1800 1600}%
\special{pa 1800 2000}%
\special{fp}%
%
\special{pn 8}%
\special{pa 1600 2000}%
\special{pa 2800 2000}%
\special{fp}%
\special{pa 1600 1600}%
\special{pa 2800 1600}%
\special{fp}%
\special{pa 2200 1600}%
\special{pa 2200 2000}%
\special{fp}%
\special{pa 2600 1600}%
\special{pa 2600 2000}%
\special{fp}%
\put(15.9000,-14.5000){\makebox(0,0){$1$}}%
\put(19.8000,-14.4000){\makebox(0,0){$1$}}%
\put(23.9000,-14.3000){\makebox(0,0){$1$}}%
\put(27.9000,-14.3000){\makebox(0,0){$1$}}%
\put(16.2000,-21.6000){\makebox(0,0){$1$}}%
\put(20.0500,-21.6500){\makebox(0,0){$1$}}%
\put(24.2500,-21.5500){\makebox(0,0){$1$}}%
\put(27.9000,-21.5000){\makebox(0,0){$1$}}%
\put(16.3000,-17.9000){\makebox(0,0){$2$}}%
\put(20.6000,-17.9000){\makebox(0,0){$2$}}%
\put(24.7000,-17.9000){\makebox(0,0){$2$}}%
\end{picture}%

\end{center}
\caption{$x$ in $(m,n)=(2,3)$.}
\end{figure}
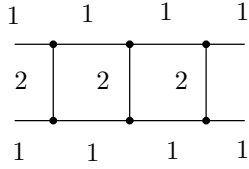

Therefore we see $L_{P^\circ}(4)\not=0$, and so that the polytope $P$ has the codegree 4. 
By Proposition \ref{bhs2}, $\dim \mathcal{P}_T(2,3)=3$. The degree of $\mathcal{P}_T(2,3)$ is 0 and $\mathcal{P}_T(2,3)$ is an unimodular simplex, and so Gorenstein.

Let $(m,n)=(2,5)$. All-one vector $\al=(1,1,\ldots,1)$ lies in $3P^\circ$. Actually, $\al$ satisfies the conditions $(1')_3$ and $(2')_3$, and easily we also confirm that, even if we take 3 or 5 points as $S$ in any choice, the condition $(3')_3$ always holds. So by Lemma \ref{rel}, $\al \in 3P^\circ$ and $tP^\circ$ is defined by $(1')_t$, $(2')_t$ and $(3')_t$. Therefore $L_{P^\circ}(3)\not=0$. Since the graph is $3$-regular, conditions $(1')_t$ and $(2')_t$ in Edmond's theorem imply that there are no lattice points in $P^\circ$ and $2P^\circ$. So the codegree of $P$ is $3$. Below we show that $L_{P^\circ}(t)=L_P(t-3)$.
 
By Lemma \ref{inj}, $\iota$ gives an injective map from $lP\cap \Z^E$ to $(l+3)P \cap \Z^E$.
We show that the inverse map $\iota^{-1}$ also gives an injective map from $(l+3)P^\circ \cap \Z^E$ to $lP\cap \Z^E$. If we could prove it, then we see that $L_{P^\circ}(l+3)=|(l+3)P^\circ \cap \Z^E|=|lP\cap \Z^E|=L_P(l)$.

Take $x\in (l+3)P\cap \Z^E$, then $y=\iota^{-1}(x)=x-\al$ satisfies $(1)_l$ and $(2)_l$.
Take $S$ so that condition $(3)_l$ does not hold, then $S$ contains no isolated points. Also by symmetry, we have only to consider $S$ with the cardinality at most half the total number of vertices. Hence, as $S$ which does not satisfy $(3)_l$, we have only four possibilities shown in Figures 2 and 3 (by considering symmetry again). In these figures, big points denote points of $S$ and thick edges denote the induced subgraph by $S$, that is, the graph consisting of the vertex set $S$ and all edges among vertices of $S$.
\begin{figure}[htbp]
\begin{center}
\unitlength 0.1in
\begin{picture}( 36.9000,  6.8000)(  9.1000,-24.0000)
%
\special{pn 20}%
\special{sh 1}%
\special{ar 1200 2000 10 10 0  6.28318530717959E+0000}%
\special{sh 1}%
\special{ar 1600 2000 10 10 0  6.28318530717959E+0000}%
\special{sh 1}%
\special{ar 2000 2000 10 10 0  6.28318530717959E+0000}%
\special{sh 1}%
\special{ar 2000 2000 10 10 0  6.28318530717959E+0000}%
%
\special{pn 8}%
\special{sh 1}%
\special{ar 910 1720 10 10 0  6.28318530717959E+0000}%
\special{sh 1}%
\special{ar 910 2290 10 10 0  6.28318530717959E+0000}%
\special{sh 1}%
\special{ar 1600 2390 10 10 0  6.28318530717959E+0000}%
\special{sh 1}%
\special{ar 2300 1720 10 10 0  6.28318530717959E+0000}%
\special{sh 1}%
\special{ar 2300 2300 10 10 0  6.28318530717959E+0000}%
\special{sh 1}%
\special{ar 2300 2300 10 10 0  6.28318530717959E+0000}%
%
\special{pn 8}%
\special{pa 910 1730}%
\special{pa 1200 2000}%
\special{fp}%
\special{pa 910 2300}%
\special{pa 1200 2000}%
\special{fp}%
\special{pa 1600 2000}%
\special{pa 1600 2400}%
\special{fp}%
\special{pa 2000 2000}%
\special{pa 2300 1730}%
\special{fp}%
\special{pa 2000 2000}%
\special{pa 2300 2300}%
\special{fp}%
%
\special{pn 20}%
\special{pa 1200 2000}%
\special{pa 2000 2000}%
\special{fp}%
%
\special{pn 8}%
\special{sh 1}%
\special{ar 3000 1800 10 10 0  6.28318530717959E+0000}%
\special{sh 1}%
\special{ar 3000 2200 10 10 0  6.28318530717959E+0000}%
%
\special{pn 20}%
\special{sh 1}%
\special{ar 3390 1790 10 10 0  6.28318530717959E+0000}%
\special{sh 1}%
\special{ar 3800 1800 10 10 0  6.28318530717959E+0000}%
\special{sh 1}%
\special{ar 4200 1800 10 10 0  6.28318530717959E+0000}%
\special{sh 1}%
\special{ar 3400 2200 10 10 0  6.28318530717959E+0000}%
\special{sh 1}%
\special{ar 3800 2200 10 10 0  6.28318530717959E+0000}%
%
\special{pn 8}%
\special{sh 1}%
\special{ar 4600 1800 10 10 0  6.28318530717959E+0000}%
\special{sh 1}%
\special{ar 4600 2200 10 10 0  6.28318530717959E+0000}%
%
\special{pn 8}%
\special{sh 1}%
\special{ar 4200 2200 10 10 0  6.28318530717959E+0000}%
\special{sh 1}%
\special{ar 4200 2200 10 10 0  6.28318530717959E+0000}%
%
\special{pn 20}%
\special{pa 3400 1800}%
\special{pa 4200 1800}%
\special{fp}%
%
\special{pn 20}%
\special{pa 3400 1800}%
\special{pa 3400 2200}%
\special{fp}%
\special{pa 3400 2200}%
\special{pa 3800 2200}%
\special{fp}%
\special{pa 3800 2200}%
\special{pa 3800 1800}%
\special{fp}%
%
\special{pn 8}%
\special{pa 3000 1800}%
\special{pa 3400 1800}%
\special{fp}%
%
\special{pn 8}%
\special{pa 3000 2200}%
\special{pa 3400 2200}%
\special{fp}%
%
\special{pn 8}%
\special{pa 3800 2200}%
\special{pa 4200 2200}%
\special{fp}%
%
\special{pn 8}%
\special{pa 4200 2200}%
\special{pa 4600 2200}%
\special{fp}%
%
\special{pn 8}%
\special{pa 4200 1800}%
\special{pa 4600 1800}%
\special{fp}%
%
\special{pn 8}%
\special{pa 4200 1800}%
\special{pa 4200 2200}%
\special{fp}%
\end{picture}%

\end{center}
\caption{$(m,n)=(2,5)$ I.}
\end{figure}
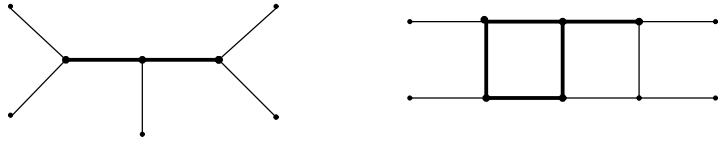
\begin{figure}[htbp]
\begin{center}
\unitlength 0.1in
\begin{picture}( 49.3300,  7.2000)( 12.0000,-28.0000)
%
\special{pn 20}%
\special{sh 1}%
\special{ar 4134 2400 10 10 0  6.28318530717959E+0000}%
\special{sh 1}%
\special{ar 4524 2400 10 10 0  6.28318530717959E+0000}%
\special{sh 1}%
\special{ar 4934 2400 10 10 0  6.28318530717959E+0000}%
\special{sh 1}%
\special{ar 5334 2400 10 10 0  6.28318530717959E+0000}%
\special{sh 1}%
\special{ar 5734 2400 10 10 0  6.28318530717959E+0000}%
%
\special{pn 8}%
\special{sh 1}%
\special{ar 3734 2400 10 10 0  6.28318530717959E+0000}%
%
\special{pn 8}%
\special{sh 1}%
\special{ar 6134 2400 10 10 0  6.28318530717959E+0000}%
\special{sh 1}%
\special{ar 6134 2400 10 10 0  6.28318530717959E+0000}%
%
\special{pn 20}%
\special{pa 3734 2400}%
\special{pa 6134 2400}%
\special{fp}%
%
\special{pn 8}%
\special{sh 1}%
\special{ar 4134 2800 10 10 0  6.28318530717959E+0000}%
\special{sh 1}%
\special{ar 4534 2800 10 10 0  6.28318530717959E+0000}%
\special{sh 1}%
\special{ar 4934 2800 10 10 0  6.28318530717959E+0000}%
\special{sh 1}%
\special{ar 5334 2800 10 10 0  6.28318530717959E+0000}%
\special{sh 1}%
\special{ar 5734 2800 10 10 0  6.28318530717959E+0000}%
%
\special{pn 8}%
\special{pa 4134 2400}%
\special{pa 4134 2800}%
\special{fp}%
%
\special{pn 8}%
\special{pa 4934 2400}%
\special{pa 4934 2800}%
\special{fp}%
%
\special{pn 8}%
\special{pa 4534 2400}%
\special{pa 4534 2800}%
\special{fp}%
%
\special{pn 8}%
\special{pa 5334 2400}%
\special{pa 5334 2800}%
\special{fp}%
\special{pa 5734 2400}%
\special{pa 5734 2800}%
\special{fp}%
\put(39.6300,-26.7000){\makebox(0,0){$a_1$}}%
\put(43.0300,-27.1000){\makebox(0,0)[lb]{$a_2$}}%
\put(47.9300,-26.7000){\makebox(0,0){$a_3$}}%
\put(50.8300,-27.2000){\makebox(0,0)[lb]{$a_4$}}%
\put(54.9300,-27.1000){\makebox(0,0)[lb]{$a_5$}}%
\put(38.8300,-22.6000){\makebox(0,0)[lb]{$b_5$}}%
\put(42.1300,-22.5000){\makebox(0,0)[lb]{$b_1$}}%
\put(47.1300,-22.0000){\makebox(0,0){$b_2$}}%
\put(50.4300,-22.5000){\makebox(0,0)[lb]{$b_3$}}%
\put(54.5300,-22.5000){\makebox(0,0)[lb]{$b_4$}}%
\put(58.8300,-22.5000){\makebox(0,0)[lb]{$b_5$}}%
%
\special{pn 20}%
\special{sh 1}%
\special{ar 1400 2400 10 10 0  6.28318530717959E+0000}%
\special{sh 1}%
\special{ar 1800 2400 10 10 0  6.28318530717959E+0000}%
\special{sh 1}%
\special{ar 2200 2400 10 10 0  6.28318530717959E+0000}%
\special{sh 1}%
\special{ar 2600 2400 10 10 0  6.28318530717959E+0000}%
\special{sh 1}%
\special{ar 3000 2400 10 10 0  6.28318530717959E+0000}%
%
\special{pn 20}%
\special{pa 1400 2400}%
\special{pa 3000 2400}%
\special{fp}%
%
\special{pn 8}%
\special{pa 1200 2200}%
\special{pa 1400 2400}%
\special{fp}%
\special{pa 1400 2400}%
\special{pa 1200 2600}%
\special{fp}%
%
\special{pn 8}%
\special{pa 1800 2400}%
\special{pa 1800 2800}%
\special{fp}%
%
\special{pn 8}%
\special{pa 2200 2400}%
\special{pa 2200 2800}%
\special{fp}%
\special{pa 2600 2400}%
\special{pa 2600 2800}%
\special{fp}%
\special{pa 3000 2400}%
\special{pa 3200 2200}%
\special{fp}%
\special{pa 3000 2400}%
\special{pa 3200 2600}%
\special{fp}%
%
\special{pn 8}%
\special{sh 1}%
\special{ar 1200 2200 10 10 0  6.28318530717959E+0000}%
\special{sh 1}%
\special{ar 1200 2600 10 10 0  6.28318530717959E+0000}%
\special{sh 1}%
\special{ar 1800 2800 10 10 0  6.28318530717959E+0000}%
\special{sh 1}%
\special{ar 2200 2800 10 10 0  6.28318530717959E+0000}%
\special{sh 1}%
\special{ar 2600 2800 10 10 0  6.28318530717959E+0000}%
\special{sh 1}%
\special{ar 3200 2600 10 10 0  6.28318530717959E+0000}%
\special{sh 1}%
\special{ar 3200 2200 10 10 0  6.28318530717959E+0000}%
\special{sh 1}%
\special{ar 3200 2200 10 10 0  6.28318530717959E+0000}%
\end{picture}%

\end{center}
\caption{$(m,n)=(2,5)$ II.}
\end{figure}
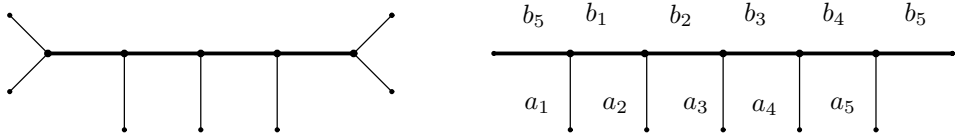

For the two cases in Figure 2 and the left case in Figure 3, we see from Corollary \ref{oyo} that the condition $(3)_l$ follows from $(1)_l$ and $(2)_l$. Next we consider the right case in Figure 3. Since $x$ satisfies $(3')_{l+3}$, $\sum_{e\in C(S,S')} x_e\ge l+4$. Therefore 
\begin{equation*}\label{x}
\sum_{e\in C(S,S')}y_e=\sum_{e\in C(S,S')}(x_e-1)\ge l-1.
\end{equation*}
If $\sum_{e\in C(S,S')}y_e<l$, then $\sum_{e\in C(S,S')}y_e=\sum_{1\le i \le 5} a_i =l-1$.
From the condition $(2)_l$ for $y$ we get
\[5l=\sum_{v\in S,\, v\in e}y_e=\sum_{1\le i \le 5} a_i+2\sum_{1\le i \le 5} b_i \equiv l-1 \pmod{2}.\]
This is a contradiction. So $y$ always satisfies the condition $(3)_l$. This implies that $\iota^{-1}$  gives an injective map from $(l+3)P^\circ \cap \Z^E$ to $lP\cap \Z^E$.

Next, in order to show the necessity, we prove the contraposition.
First let $m=2$ and $n\ge 7$ be odd. Similar to the case when $(m,n)=(2,5)$, we see easily that $\al$ lies in $3P^\circ$ and $tP^\circ$ is defined by $(1')_t$, $(2')_t$ and $(3')_t$. So the codegree is $3$. It is also similar that $\iota$ gives an injective map from $lP\cap \Z^E$ to $(l+3)P^\circ \cap \Z^E$. If we could prove the existence of $y$ such that $y \not\in 2P$ and $\iota(y)\in 5P^\circ$, then we see $L_P(2)<L_{P^\circ}(5)$, and so the polytope is not Gorenstein. Define a vector $y$ as in Figure 4. 
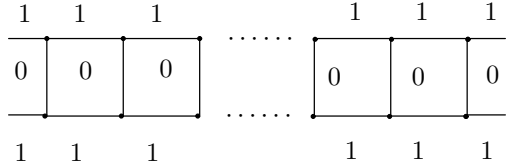
\begin{figure}[htbp]
\begin{center}
\unitlength 0.1in
\begin{picture}( 26.0000,  7.4000)( 32.0000,-18.7000)
\put(49.8000,-13.0000){\makebox(0,0)[lb]{$1$}}%
%
\special{pn 8}%
\special{sh 1}%
\special{ar 3400 1400 10 10 0  6.28318530717959E+0000}%
\special{sh 1}%
\special{ar 3800 1400 10 10 0  6.28318530717959E+0000}%
\special{sh 1}%
\special{ar 4200 1400 10 10 0  6.28318530717959E+0000}%
\put(45.3000,-14.0000){\makebox(0,0){$\cdots \cdots$}}%
%
\special{pn 8}%
\special{sh 1}%
\special{ar 4804 1404 10 10 0  6.28318530717959E+0000}%
%
\special{pn 8}%
\special{sh 1}%
\special{ar 5204 1404 10 10 0  6.28318530717959E+0000}%
%
\special{pn 8}%
\special{sh 1}%
\special{ar 5604 1404 10 10 0  6.28318530717959E+0000}%
\special{sh 1}%
\special{ar 5604 1404 10 10 0  6.28318530717959E+0000}%
%
\special{pn 8}%
\special{pa 3400 1400}%
\special{pa 3800 1400}%
\special{fp}%
\special{pa 4200 1400}%
\special{pa 3800 1400}%
\special{fp}%
%
\special{pn 8}%
\special{pa 4800 1400}%
\special{pa 5200 1400}%
\special{fp}%
\special{pa 5200 1400}%
\special{pa 5600 1400}%
\special{fp}%
%
\special{pn 8}%
\special{sh 1}%
\special{ar 3390 1806 10 10 0  6.28318530717959E+0000}%
\special{sh 1}%
\special{ar 3790 1806 10 10 0  6.28318530717959E+0000}%
\special{sh 1}%
\special{ar 4190 1806 10 10 0  6.28318530717959E+0000}%
\put(45.2000,-18.0500){\makebox(0,0){$\cdots \cdots$}}%
%
\special{pn 8}%
\special{sh 1}%
\special{ar 4794 1808 10 10 0  6.28318530717959E+0000}%
%
\special{pn 8}%
\special{sh 1}%
\special{ar 5194 1808 10 10 0  6.28318530717959E+0000}%
%
\special{pn 8}%
\special{sh 1}%
\special{ar 5594 1808 10 10 0  6.28318530717959E+0000}%
\special{sh 1}%
\special{ar 5594 1808 10 10 0  6.28318530717959E+0000}%
%
\special{pn 8}%
\special{pa 3390 1806}%
\special{pa 3790 1806}%
\special{fp}%
\special{pa 4190 1806}%
\special{pa 3790 1806}%
\special{fp}%
%
\special{pn 8}%
\special{pa 4790 1806}%
\special{pa 5190 1806}%
\special{fp}%
\special{pa 5190 1806}%
\special{pa 5590 1806}%
\special{fp}%
%
\special{pn 8}%
\special{pa 3400 1400}%
\special{pa 3400 1800}%
\special{fp}%
%
\special{pn 8}%
\special{pa 3800 1400}%
\special{pa 3800 1800}%
\special{fp}%
%
\special{pn 8}%
\special{pa 4200 1400}%
\special{pa 4200 1800}%
\special{fp}%
%
\special{pn 8}%
\special{pa 5600 1400}%
\special{pa 5600 1800}%
\special{fp}%
%
\special{pn 8}%
\special{pa 5200 1400}%
\special{pa 5200 1800}%
\special{fp}%
%
\special{pn 8}%
\special{pa 4800 1400}%
\special{pa 4800 1800}%
\special{fp}%
%
\special{pn 8}%
\special{pa 3200 1400}%
\special{pa 3400 1400}%
\special{fp}%
%
\special{pn 8}%
\special{pa 3400 1800}%
\special{pa 3200 1800}%
\special{fp}%
%
\special{pn 8}%
\special{pa 5600 1400}%
\special{pa 5800 1400}%
\special{fp}%
%
\special{pn 8}%
\special{pa 5600 1800}%
\special{pa 5800 1800}%
\special{fp}%
\put(32.5000,-13.1000){\makebox(0,0)[lb]{$1$}}%
\put(39.4000,-13.1000){\makebox(0,0)[lb]{$1$}}%
\put(35.4000,-13.1000){\makebox(0,0)[lb]{$1$}}%
\put(53.3000,-13.0000){\makebox(0,0)[lb]{$1$}}%
\put(56.9000,-13.0000){\makebox(0,0)[lb]{$1$}}%
\put(32.3000,-16.1000){\makebox(0,0)[lb]{$0$}}%
\put(35.7000,-16.1000){\makebox(0,0)[lb]{$0$}}%
\put(39.9000,-16.0000){\makebox(0,0)[lb]{$0$}}%
\put(48.7000,-16.4000){\makebox(0,0)[lb]{$0$}}%
\put(53.1000,-16.4000){\makebox(0,0)[lb]{$0$}}%
\put(57.0000,-16.3000){\makebox(0,0)[lb]{$0$}}%
\put(49.6000,-20.3000){\makebox(0,0)[lb]{$1$}}%
\put(32.3000,-20.4000){\makebox(0,0)[lb]{$1$}}%
\put(39.2000,-20.4000){\makebox(0,0)[lb]{$1$}}%
\put(35.2000,-20.4000){\makebox(0,0)[lb]{$1$}}%
\put(53.1000,-20.3000){\makebox(0,0)[lb]{$1$}}%
\put(56.7000,-20.3000){\makebox(0,0)[lb]{$1$}}%
\end{picture}%

\end{center}
\caption{$y$ for $(m,n)=(2,n),\,n\ge 7$.}
\end{figure}

When we take all points on the upper row as $S$, $y$ does not satisfy the condition $(3)_2$.  So $y\not\in 2P$. We show that $x=\iota(y)=y+\al\in 5P^\circ$. If $\sum_{e\in C(S,S')}y_e \ge t$ for $S$, then since $\al$ satisfies $(3')_3$, $\sum_{e\in C(S,S')}x_e > t+3$. Thus, we have only to consider $S$ such that $|S|\le n$ is odd and that $\sum_{e\in C(S,S')}y_e < 2$. 
The inequality $\sum_{e\in C(S,S')}y_e < 2$ holds only when we take all points on the upper row as $S$. 
In this case, since $n$ edges goes from the upper row to the lower row, so $n=\sum_{e\in C(S,S')}x_e \ge t+4=6$. Therefore $x\in 5P^\circ$.

Next we consider the case when $m\ge 4$ is even and $n=3$. By symmetry we put $m=3$ and let $n\ge 4$ be even. Since the graph is $4$-regular, $\al$ satisfies $(1')_4$ and $(2')_4$. Also from Lemma \ref{first} we see that $\al$ satisfies the condition $(3')_4$, too. Therefore by Lemma \ref{rel}, $\al \in 4P^\circ$ and $tP^\circ$ is defined by $(1')_t$, $(2')_t$ and $(3')_t$, and so the codegree of $P$ is $4$. By Lemma \ref{inj},  $\iota$ gives an injective map from $lP\cap \Z^E$ to $(l+4)P^\circ \cap \Z^E$. Therefore, in order to prove that the polytope is not Gorenstein, it is sufficient to prove the existence of $y$ such that $y\not\in 3P$ and $\iota(y)\in 7P^\circ$. Define a vector $y$ as in Figure 5. Here we assign 0 to edges except for the thick ones. 

\begin{figure}[htbp]
\begin{center}
\unitlength 0.1in
\begin{picture}( 48.0000, 12.7500)( 26.0000,-30.0000)
%
\special{pn 8}%
\special{sh 1}%
\special{ar 2800 2000 10 10 0  6.28318530717959E+0000}%
\special{sh 1}%
\special{ar 3200 2000 10 10 0  6.28318530717959E+0000}%
\special{sh 1}%
\special{ar 3600 2000 10 10 0  6.28318530717959E+0000}%
\special{sh 1}%
\special{ar 4000 2000 10 10 0  6.28318530717959E+0000}%
\special{sh 1}%
\special{ar 4400 2000 10 10 0  6.28318530717959E+0000}%
\special{sh 1}%
\special{ar 4800 2000 10 10 0  6.28318530717959E+0000}%
%
\special{pn 8}%
\special{pa 2600 2000}%
\special{pa 5000 2000}%
\special{fp}%
%
\special{pn 8}%
\special{sh 1}%
\special{ar 2800 2400 10 10 0  6.28318530717959E+0000}%
\special{sh 1}%
\special{ar 3200 2400 10 10 0  6.28318530717959E+0000}%
\special{sh 1}%
\special{ar 3600 2400 10 10 0  6.28318530717959E+0000}%
\special{sh 1}%
\special{ar 4000 2400 10 10 0  6.28318530717959E+0000}%
\special{sh 1}%
\special{ar 4400 2400 10 10 0  6.28318530717959E+0000}%
\special{sh 1}%
\special{ar 4800 2400 10 10 0  6.28318530717959E+0000}%
%
\special{pn 8}%
\special{pa 2600 2400}%
\special{pa 5000 2400}%
\special{fp}%
%
\special{pn 8}%
\special{sh 1}%
\special{ar 2800 2800 10 10 0  6.28318530717959E+0000}%
\special{sh 1}%
\special{ar 3200 2800 10 10 0  6.28318530717959E+0000}%
\special{sh 1}%
\special{ar 3600 2800 10 10 0  6.28318530717959E+0000}%
\special{sh 1}%
\special{ar 4000 2800 10 10 0  6.28318530717959E+0000}%
\special{sh 1}%
\special{ar 4400 2800 10 10 0  6.28318530717959E+0000}%
\special{sh 1}%
\special{ar 4800 2800 10 10 0  6.28318530717959E+0000}%
%
\special{pn 8}%
\special{pa 2600 2800}%
\special{pa 5000 2800}%
\special{fp}%
%
\special{pn 8}%
\special{pa 2800 1800}%
\special{pa 2800 3000}%
\special{fp}%
%
\special{pn 8}%
\special{pa 3200 1800}%
\special{pa 3200 3000}%
\special{fp}%
%
\special{pn 8}%
\special{pa 3600 1800}%
\special{pa 3600 3000}%
\special{fp}%
%
\special{pn 8}%
\special{pa 3600 1800}%
\special{pa 3600 3000}%
\special{fp}%
%
\special{pn 8}%
\special{pa 4000 1800}%
\special{pa 4000 3000}%
\special{fp}%
%
\special{pn 8}%
\special{pa 4400 1800}%
\special{pa 4400 3000}%
\special{fp}%
%
\special{pn 8}%
\special{pa 4400 1800}%
\special{pa 4400 3000}%
\special{fp}%
%
\special{pn 8}%
\special{pa 4800 1800}%
\special{pa 4800 3000}%
\special{fp}%
%
\special{pn 8}%
\special{pa 4800 1800}%
\special{pa 4800 3000}%
\special{fp}%
%
\special{pn 20}%
\special{pa 4400 2000}%
\special{pa 4400 2400}%
\special{fp}%
\special{pa 4400 2400}%
\special{pa 4800 2400}%
\special{fp}%
\special{pa 4800 2400}%
\special{pa 4800 2000}%
\special{fp}%
\special{pa 4800 2000}%
\special{pa 4400 2000}%
\special{fp}%
%
\special{pn 20}%
\special{pa 4800 1800}%
\special{pa 4800 3000}%
\special{fp}%
%
\special{pn 20}%
\special{pa 3600 2000}%
\special{pa 4000 2000}%
\special{fp}%
%
\special{pn 20}%
\special{pa 3600 2400}%
\special{pa 4000 2400}%
\special{fp}%
%
\special{pn 20}%
\special{pa 2800 2000}%
\special{pa 3200 2000}%
\special{fp}%
\special{pa 2800 2400}%
\special{pa 3200 2400}%
\special{fp}%
%
\special{pn 20}%
\special{pa 4000 2800}%
\special{pa 4400 2800}%
\special{fp}%
\special{pa 3200 2800}%
\special{pa 3600 2800}%
\special{fp}%
\special{pa 2800 2800}%
\special{pa 2600 2800}%
\special{fp}%
%
\special{pn 8}%
\special{sh 1}%
\special{ar 7200 2000 10 10 0  6.28318530717959E+0000}%
\special{sh 1}%
\special{ar 6800 2000 10 10 0  6.28318530717959E+0000}%
\special{sh 1}%
\special{ar 6400 2000 10 10 0  6.28318530717959E+0000}%
\special{sh 1}%
\special{ar 6000 2000 10 10 0  6.28318530717959E+0000}%
\special{sh 1}%
\special{ar 5600 2000 10 10 0  6.28318530717959E+0000}%
\special{sh 1}%
\special{ar 5200 2000 10 10 0  6.28318530717959E+0000}%
%
\special{pn 8}%
\special{pa 7400 2000}%
\special{pa 5000 2000}%
\special{fp}%
%
\special{pn 8}%
\special{sh 1}%
\special{ar 7200 2400 10 10 0  6.28318530717959E+0000}%
\special{sh 1}%
\special{ar 6800 2400 10 10 0  6.28318530717959E+0000}%
\special{sh 1}%
\special{ar 6400 2400 10 10 0  6.28318530717959E+0000}%
\special{sh 1}%
\special{ar 6000 2400 10 10 0  6.28318530717959E+0000}%
\special{sh 1}%
\special{ar 5600 2400 10 10 0  6.28318530717959E+0000}%
\special{sh 1}%
\special{ar 5200 2400 10 10 0  6.28318530717959E+0000}%
%
\special{pn 8}%
\special{pa 7400 2400}%
\special{pa 5000 2400}%
\special{fp}%
%
\special{pn 8}%
\special{sh 1}%
\special{ar 7200 2800 10 10 0  6.28318530717959E+0000}%
\special{sh 1}%
\special{ar 6800 2800 10 10 0  6.28318530717959E+0000}%
\special{sh 1}%
\special{ar 6400 2800 10 10 0  6.28318530717959E+0000}%
\special{sh 1}%
\special{ar 6000 2800 10 10 0  6.28318530717959E+0000}%
\special{sh 1}%
\special{ar 5600 2800 10 10 0  6.28318530717959E+0000}%
\special{sh 1}%
\special{ar 5200 2800 10 10 0  6.28318530717959E+0000}%
%
\special{pn 8}%
\special{pa 7400 2800}%
\special{pa 5000 2800}%
\special{fp}%
%
\special{pn 8}%
\special{pa 7200 1800}%
\special{pa 7200 3000}%
\special{fp}%
%
\special{pn 8}%
\special{pa 6800 1800}%
\special{pa 6800 3000}%
\special{fp}%
%
\special{pn 8}%
\special{pa 6400 1800}%
\special{pa 6400 3000}%
\special{fp}%
%
\special{pn 8}%
\special{pa 6400 1800}%
\special{pa 6400 3000}%
\special{fp}%
%
\special{pn 8}%
\special{pa 6000 1800}%
\special{pa 6000 3000}%
\special{fp}%
%
\special{pn 8}%
\special{pa 5600 1800}%
\special{pa 5600 3000}%
\special{fp}%
%
\special{pn 8}%
\special{pa 5600 1800}%
\special{pa 5600 3000}%
\special{fp}%
%
\special{pn 8}%
\special{pa 5200 1800}%
\special{pa 5200 3000}%
\special{fp}%
%
\special{pn 8}%
\special{pa 5200 1800}%
\special{pa 5200 3000}%
\special{fp}%
%
\special{pn 20}%
\special{pa 5600 2000}%
\special{pa 5600 2400}%
\special{fp}%
\special{pa 5600 2400}%
\special{pa 5200 2400}%
\special{fp}%
\special{pa 5200 2400}%
\special{pa 5200 2000}%
\special{fp}%
\special{pa 5200 2000}%
\special{pa 5600 2000}%
\special{fp}%
%
\special{pn 20}%
\special{pa 5200 1800}%
\special{pa 5200 3000}%
\special{fp}%
%
\special{pn 20}%
\special{pa 6400 2000}%
\special{pa 6000 2000}%
\special{fp}%
%
\special{pn 20}%
\special{pa 6400 2400}%
\special{pa 6000 2400}%
\special{fp}%
%
\special{pn 20}%
\special{pa 7200 2000}%
\special{pa 6800 2000}%
\special{fp}%
\special{pa 7200 2400}%
\special{pa 6800 2400}%
\special{fp}%
%
\special{pn 20}%
\special{pa 6000 2800}%
\special{pa 5600 2800}%
\special{fp}%
\special{pa 6800 2800}%
\special{pa 6400 2800}%
\special{fp}%
\special{pa 7200 2800}%
\special{pa 7400 2800}%
\special{fp}%
%
\special{pn 8}%
\special{pa 4800 2800}%
\special{pa 5210 2790}%
\special{fp}%
%
\special{pn 20}%
\special{pa 4800 2800}%
\special{pa 5200 2800}%
\special{fp}%
\put(42.7000,-22.0000){\makebox(0,0){2}}%
\put(45.8000,-18.7000){\makebox(0,0){1}}%
\put(45.8000,-22.9000){\makebox(0,0){1}}%
\put(48.9000,-22.4000){\makebox(0,0){1}}%
\put(49.0000,-26.2000){\makebox(0,0){1}}%
\put(48.9000,-29.9000){\makebox(0,0){1}}%
\put(50.6000,-18.2000){\makebox(0,0){1}}%
\put(51.0000,-22.3000){\makebox(0,0){1}}%
\put(51.0000,-26.2000){\makebox(0,0){1}}%
\put(50.9000,-29.8000){\makebox(0,0){1}}%
\put(57.2000,-22.0000){\makebox(0,0){2}}%
\put(54.1500,-18.6500){\makebox(0,0){1}}%
\put(53.8000,-22.8000){\makebox(0,0){1}}%
\put(29.7000,-18.7000){\makebox(0,0){3}}%
\put(37.9000,-18.6000){\makebox(0,0){3}}%
\put(29.7000,-22.7000){\makebox(0,0){3}}%
\put(37.9000,-22.6000){\makebox(0,0){3}}%
\put(34.1000,-26.4000){\makebox(0,0){3}}%
\put(41.8000,-26.4000){\makebox(0,0){3}}%
\put(26.5000,-26.6000){\makebox(0,0){3}}%
\put(58.3000,-26.5000){\makebox(0,0){3}}%
\put(66.0000,-26.7000){\makebox(0,0){3}}%
\put(62.0000,-22.4000){\makebox(0,0){3}}%
\put(62.0000,-18.1000){\makebox(0,0){3}}%
\put(69.7000,-18.2000){\makebox(0,0){3}}%
\put(69.9000,-22.3000){\makebox(0,0){3}}%
\put(73.6000,-26.5000){\makebox(0,0){3}}%
\put(49.9000,-29.3000){\makebox(0,0){1}}%
\put(43.0000,-18.8000){\makebox(0,0){$c_1$}}%
\put(42.8000,-24.9000){\makebox(0,0){$c_2$}}%
\put(47.1000,-25.2000){\makebox(0,0){$c_4$}}%
\put(46.7000,-29.2000){\makebox(0,0){$c_5$}}%
\put(46.9000,-18.5000){\makebox(0,0){$c_3$}}%
\put(57.1000,-18.9000){\makebox(0,0){$d_1$}}%
\put(57.1000,-25.0000){\makebox(0,0){$d_2$}}%
\put(52.9500,-18.7500){\makebox(0,0){$d_3$}}%
\put(53.0000,-25.1000){\makebox(0,0){$d_4$}}%
\put(52.9000,-29.0000){\makebox(0,0){$d_5$}}%
\put(48.8000,-18.2000){\makebox(0,0){1}}%
\end{picture}%

\end{center}
\caption{$y$ for $(m,n)=(3,n),\,n\ge 4$.}
\end{figure}

Take $c_i$ 's $(1\le i \le 5)$ as $S$, then $\sum_{e\in C(S,S')}y_e=1<3$. So $y\not\in 3P$. Next, we show that $x=\iota(y)\in 7P^\circ$. It is clear that $(1')_7$ and $(2')_7$ hold. We have to show that $(3')_7$ holds for any $S$ with the odd cardinality. In a similar way to the case when $m=2,\,n\ge 7$, using Lemma \ref{first}, we see that it is sufficient to consider only $S$ such that $\sum_{e\in C(S,S')}y_e \le 1$. To choose points of $S$ satisfying $\sum_{e\in C(S,S')}y_e \le 1$ we need to take all or none of $c_i$'s. This is also similar for ${d_i}$'s. Therefore, candidates for $S$ consist of all $c_i$'s, none of $d_i$'s and some matchings assigned by 3. In this case, we show that $|C(S,S')|\ge 8$.
If we could show this fact, we see that $\sum_{e\in C(S,S')}x_e=\sum_{e\in C(S,S')}(y_e+1)\ge 1+\sum_{e\in C(S,S')}1\ge 9$, and so $x \in 7P^\circ$. In this case, we see that each row has at least 2 bridges, and since the matching in the third row is out of sync with the first and second ones, bridges traverse from the third row to the first and second rows. So $|C(S,S')|\ge 8$.

Next we consider the case when $m \ge 4$ is even  and $n\ge 5$ is odd. From Lemmas \ref{first} and \ref{rel}, we see that $\al\in 4P^\circ$ and $tP^\circ$ is defined by $(1')_t$, $(2')_t$ and $(3')_t$, and so that the codegree of $P$ is 4. By Lemma \ref{inj}, $\iota$ gives an injective map from $lP\cap \Z^E$ to $(l+4)P^\circ \cap \Z^E$.  So we also construct a vector $y$ such that $y\not\in 2P$ and $x=\iota(y) \in 6P^\circ$. Define a vector $y$ so that each horizontal edge gets a 0, and each vertical edge gets 0. Take all points on the first row as $S$, $\sum_{e\in C(S,S')} y_e=0<2$. So $y \not\in 2P$. 

We show that $x=\iota(y)$ satisfies $\sum_{e\in C(S,S')} x_e >6$ for any $S$ with the cardinality odd. By Lemma \ref{first} and similar arguments as the above, we consider as candidates for $S$ only ones such that $\sum_{e\in C(S,S')} y_e=0$. To choose $S$ satisfying this condition, we need to take all points or none of points on each row. Therefore in this case, there are $2n$ bridges. So we have $|C(S,S')|\ge 2n\ge 10$, and so $\sum_{e \in C(S,S')} x_e=\sum_{e \in C(S,S')} (y_e+1)\ge |C(S,S')|\ge 10$. So we get $x\in 6P^\circ$.   
\end{proof}
\begin{remark}
By Theorem \ref{hokan}, $P=\mathcal{P}_T(2,5)$ is a Gorenstein polytope of codegree $3$. We also know $\dim \mathcal{P}_T(2,5)=5$ by Proposition \ref{bhs2}, and so the degree of $P$ is $3$. The vertices of $\mathcal{P}_T(2,5)$ coincide with lattice points in $\mathcal{P}_T(2,5)$, and they are given by perfect matchings of $\mathcal{G}_T(2,5)$. We can count easily the number of perfect matchings of $\mathcal{G}_T(2,5)$ which is equal to 11. Since $h_1=L_P(1)-(d+1)=11-6=5$, the Ehrhart series is 
\[\Ehr_P(z)=\frac{1+5z+5z^2+z^3}{(1-z)^6}.\] 
\end{remark}
\section{S-matching polytope}\label{smp}
In this section, we construct new polytopes from graphs. 

Let $G=(V,E)$ be a graph, and for a subset $S\subset V$, we denote by $\langle S\rangle$ the induced subgraph by $S$ in $G$, that is, the graph consisting of the vertex set $S$ and all edges among vertices of $S$. We also define a subgraph $N_G(S)=(V_S,E_S)$ of $G$ as follows:
\[\Gamma(S):=\{x\in S'\mid \mbox{$(x,y)\in E$ for some $y \in S$}\},\]
\[V_S:=S\cup \Gamma(S),\;E_S:=C(S,V).\]
We call $N_G(S)$ the \textit{neighbor graph} of $S$ in $G$. We call $M\subset E_S$ an $S$-\textit{matching} if any two distinct edges in $M$ do not meet at points in $S$, and any point of $S$ lies on some edge in $M$. 
\begin{defi}[S-matching polytope]
Let $N_G(S)=(V_S,E_S)$ be the neighbor graph of a subset $S$. Then we define the \textit{S-matching polytope} $P_S$ of $S$ to be the convex hull in $\R^{E_S}$ of the characteristic vectors of all $S$-matchings:
\[P_S:=\con\{\chi_M\in \R^{E_S}\mid \mbox{$M$ is an $S$-matching}\},\]
where $\chi_M \in \R^{E_S}$ denotes the characteristic vector of $M$ as defined in Section \ref{intro}. 
\end{defi}
\begin{remark}
When we take all points of $V$ as $S$, $P_S$ coincides with the perfect matching polytope of $G$.
\end{remark}
\begin{theo}\label{main}
Let $G=(V,E)$ be a graph, and for a subset $S\subset V$, let $N_G(S)=(V_S,E_S)$ be the neighbor graph. Assume that $\langle S\rangle$ is bipartite. Then $x\in \R^{E_S}$ is in $P_S$ if and only if
 the following conditions hold:
 \begin{enumerate}[$(1)$]
\item $x_e\ge 0 \;(\forall e \in E_S)$,
\item $\sum_{v \in e}x_e=1$ $(\forall v\in S)$.
\end{enumerate}
\end{theo}
\begin{proof}
The proof is similar to Vempala \cite{vempala}'s proof in his lecture note on Edmond's theorem for bipartite graphs. 

Denote by $C$ the polytope defined by the above conditions (1) and (2). Let $M$ be an $S$-matching, then clearly $\chi_M$ satisfies (1) and (2). Therefore $P_S\subset C$. Conversely if we take a lattice point in $C$, then there is an $S$-matching corresponding to the lattice point. So it is sufficient to show that all vertices of $C$ are lattice points.

Assume that there exists a non-integral vertex $x=(\ldots,x_e,\ldots)$ of $C$. For such a $x$, we define a subgraph $N_G(S)_x$ of $N_G(S)$ as follows: the vertex set of $N_G(S)_x$ is $V_S$, and the edge set consists of all $e$'s in $E_S$ such that $x_e$ is not integral. We divide the cases into  (I) the case when $N_G(S)_x$ does not contain cycles, and  (II) the case when $N_G(S)_x$ contains cycles.

(I) In this case, since a connected component of $N_G(S)_x$ is a tree, there are at least two points of degree 1 in a connected component. The points of degree $1$ cannot be in $S$ by condition (2). Let $e_1,\, e_2,\ldots,e_m$ be a path connecting two such points of degree 1, and let us define $\epsilon:=\min\{x_{e_i}, 1-x_{e_i}\mid 1\le i \le m\}$. Then we define $\underbar{x}$ as follows: for $e$ not in the path, let $\underbar{x}_{e}=x_e$, and $\underbar{x}_{e_1}:=x_{e_1}-\epsilon$, $\underbar{x}_{e_2}:=x_{e_2}+\epsilon, \ldots$, that is, we alternately add and subtract $\epsilon$ to the $x_{e_i}$'s in starting from the subtraction. Also define $\overline{x}$ in a similar way except starting from the addition by $\epsilon$. Then clearly both $\underbar{x}$ and $\overline{x}$ still satisfy conditions (1) and (2).
Therefore $\underbar{x},\,\overline{x}\in C$. Since $x=(\underbar{x}+\overline{x})/2$, $x$ cannot be a vertex of $C$. This is a contradiction.

(II) In this case, there are cycles in $N_G(S)_x$. If there exists a cycle of even length, then in  a similar way to the case (I), we can define $\overline{x}$ and $\underbar{x}$ in $C$, and $x$ cannot be a vertex of $C$. 

If there exists a cycle of odd length, then since $\langle S\rangle$ is bipartite, the cycle must contain a point of $\Gamma(S)$. Let the cycle be $e_1, \ldots, e_m$ and $v\in \Gamma(S)$ be incident to $e_1$ and $e_m$. Then we define $\overline{x}$ as follows: for $e$ not in the cycle, $\overline{x}_e=x_e$ and $\overline{x}_{e_1}=x_{e_1}+\epsilon$, $\overline{x}_{e_2}=x_{e_2}-\epsilon,\ldots$, that is, we alternately add and subtract $\epsilon$ to $x_{e_i}$'s in starting from the addition. This is the same to the case (I). But finally we have $\overline{x}_{e_m}=x_{e_m}+\epsilon$. Since the common point between $e_1$ and $e_m$ is $v\in \Gamma(S)$, the addtion by $\epsilon$ does not violate condition (2). So $\overline{x}\in C$. We can define $\underbar{x}\in C$ in a similar way to $\overline{x}$ except starting from acting the subtraction. Since $\underbar{x},\,\overline{x}\in C$ and $x=(\underbar{x}+\overline{x})/2$, $x$ cannot be a vertex of $C$ in this case, and again we get a contradiction.
\end{proof}

Next we give a formula for the dimension of $P_S$ in the case when $\langle S\rangle$ is bipartite.
\begin{proposition}\label{dim}
Let $G=(V,E)$ be a graph, $S \subset V$, $N_G(S)=(V_S,E_S)$ be the neighbor graph. 
Assume that $\langle S\rangle$ is connected and bipartite, and that any $e\in E_S$ lies in some $S$-matching. Then
\[\dim P_S=
\begin{cases} 
|E_S|-|S|: & \mbox{if $\Gamma(S)\not=\emptyset$},\\
|E_S|-|S|+1: & \mbox{otherwise}.
\end{cases}\]
\end{proposition}
\begin{proof}
Define an $S\times E_S$ matrix $R$ as follows: for $v\in S$ and $e\in E_S$,
\[R_{v,e}=
\begin{cases}
1: & \mbox{if $v\in e$},\\
0: & \mbox{otherwise}.
\end{cases}\]
By Theorem \ref{main}, $P_S$ coincides with the solution space of $Rx={^t\al}$ and $x\in \R_{\ge 0}^{E_S}$, where $^tx$ is the transpose of $x$. For a suitable $\{\lambda_M\}$ such that $\sum \lambda_M=1,\,\lambda_M >0$, we put $x=\sum \lambda_M\chi_M$ where the sum runs through all $S$-matchings. By the assumption that any $e\in E_S$ lies in some $S$-matching, we have a solution $x$ of $Rx={^t\al}$ with $x_e>0 \;(\forall e\in E_S)$. So the dimension of $P_S$ is equal to the dimension of the solution space for $Rx=0$. Below we investigate the rank of $R$. Denote the $v$-th row vector and the $e$-th column vector of $R$ by $R_v$ and $R^{(e)}$, respectively. The columns of $R$ are divided into ones of $C(S,S)$ and $C(S,S')$.

First consider the case when $\Gamma(S)=\emptyset$. Then there are no columns in $C(S,S')$. 
Assume that $\sum_{v} a_vR_v=0$. If $e=(v,v')\in E_S$, then $R^{(e)}$ has the entries with 1 at $v$ and $v'$, and with 0 at the other points of $S$. Therefore, since $\sum_{v} a_vR_v=0$, we have $a_v=-a_{v'}$. So we see that if $(v,v')\in E_S$, then $a_v=-a_{v'}$. Since $\langle S\rangle$ is bipartite, there are disjoint subsets $S_1$ and $S_2$ such that $S=S_1\cup S_2$ and $C(S_i,S_i)=\emptyset$. Since $\langle S\rangle$ is connected, if $a_v=\lambda$ for some $v\in S_1$, 
then $a_u=\lambda$ for any $u\in S_1$ and $a_w=-\lambda$ for any $w\in S_2$. This implies that the dimension of the row space of $R$ is $|S|-1$. So the dimension of the solution space of $Rx=0$ is equal to $|E_S|-|S|+1$.

Next we consider the case when $\Gamma(S)\not=\emptyset$. Assume $\sum_{v} a_vR_v=0$. By the definition of $\Gamma(S)$, for $v'\in \Gamma(S)$, there exists some $v\in S$ such that $e=(v,v')\in E_S$. Then the entries of $R^{(e)}$ is 1 at $v$ and 0 at the other positions. Therefore $a_v=0$. While, it is the same to the above case that if $v,\,v'\in S$ are adjacent, then $a_v=-a_v'$. Since $\langle S\rangle$ is connected, we get, if $\sum_{v} a_vR_v=0$, then $a_v=0\;(\forall v \in S)$, and so, the dimension of the row space of $R$ is $|S|$. Therefore the dimension of the solution space for $Rx=0$ is $|E_S|-|S|$.
\end{proof}
\begin{remark}
Let us consider the case when $\langle S \rangle$ is unconnected. Then denote the connected component of $\langle S \rangle$ by $C_1,\ldots,C_k$ and put each S-matching polytope to be $P_{C_i}$. Then we have $P_S=P_{C_1}\times P_{C_2}\times\cdots\times P_{C_k}$. So the dimension is
\[\dim P_S=\sum_{i=1}^k \dim P_{C_i}.\]
\end{remark}
\begin{cor}\label{gorcor}
Assume that the induced subgraph $\langle S\rangle$ is bipartite, and that any $v\in S$ has a constant degree $k$. Then $P_S$ is a Gorenstein polytope of codegree $k$.
\end{cor}
\begin{proof}
The proof is similar to the one of Beck et. al. \cite{BHS}. By Theorem \ref{main},  $x\in tP_S^\circ$ if and only if the following conditions hold:
\begin{enumerate}[(1")]
\item $x_e > 0\;(\forall e \in E_S)$, 
\item $\sum_{v\in e} x_e=t\;(\forall v \in S)$. 
\end{enumerate}
By conditions (1") and (2"), there are no lattice points in $tP_S^\circ$ for $t<k$, and $\al\in\R^{E_S}$ is a unique lattice point of $kP_S^\circ$. So $L_{P_S}(k)=1$. 

Consider $\iota:\R^{E_S}\longrightarrow \R^{E_S}$: $\iota(y)=y+\al$. Since the degree for any point of $S$ is a constant $k$, the number of variables appearing in the conditions (2) and (2") is always $k$. Therefore for $y\in tP_S\cap \Z^{E_S}$, 
\[\sum_{v\in e} \iota(y)_e=\sum_{v\in e}(y_e+1)=t+k,\]
 and so $\iota(x)\in (t+k)P_S\cap \Z^{E_S}$. Conversely, let $x\in (t+k)P_S^\circ\cap \Z^{E_S}$. Since $x_e\ge 1$, we get $\iota^{-1}(x)_e=x_e-1\ge 0$. Since
\[\sum_{v\in e} \iota^{-1}(x)_e=\sum_{v\in e}(x_e-1)=t+k-k=t,\]
we also have $\iota^{-1}(x)\in tP_S\cap \Z^{E_S}$. 

After all, $\iota$ gives an one-to-one correspondence between $tP_S\cap \Z^{E_S}$ and $(t+k)P_S\cap \Z^{E_S}$. Therefore $L_{P_S}(t)=L_{P_S^\circ}(t+k)$, and so  $P_S$ is Gorenstein.
\end{proof}
\begin{remark}
As mentioned in Beck et. al. \cite{BHS}, combining the results of Ohsugi-Hibi \cite{oh}, Sullivant \cite{sull}, Athanasiadis \cite{Athan} and Bruns-R\"{o}mer \cite{BR}, we see that the S-matching polytopes in Corollary \ref{gorcor} are compressed  Gorenstein polytopes, and so the coefficients of the $h^*$-polynomials are unimodal.
\end{remark}
\begin{example}
By Beck et. al. \cite{BHS}, we know that for sufficient large $m$ and $n$ , $\mathcal{P}(m,n)$ is not Gorenstein (or even for the case when $mn$ is odd since in that case there are no perfect matchings). For the $m\times n$ grid graph $G=\mathcal{G}(m,n)=(V,E)$, put
\[S:=\{(i,j)\mid 1\le i\le m-2,\, 1\le j\le n-2\}\subset V,\]
then clearly $\langle S\rangle$ is bipartite and any point of $S$ has the degree $4$.
Thus, by Corollary \ref{gorcor} the S-matching polytope $P_S$ is a Gorenstein polytope of codegree $4$.
Also by Proposition \ref{dim} $\dim P_S=|E_S|-|S|=mn-m-n$.
\end{example}
\begin{cor}\label{oyo}
Let $G=(V,E)$ be a graph, and $S \subset V$ with $|S|$ odd and $\langle S \rangle$ bipartite. Let $t\in \R_{\ge 0}$ and $x\in\R_{\ge 0}^{E}$ such that $\sum_{v\in e} x_e=t\;(\forall v\in S)$. Then $\sum_{e\in C(S,S')} x_e\ge t$.
\end{cor}
\begin{proof}
By considering a projection we may assume that $x\in\R_{\ge 0}^{E_S}$. Since $\sum_{v\in e} x_e=t\;(\forall v\in S)$, by Theorem \ref{main}, $x\in tP_S$. The vertices of $P_S$ correspond  to the characteristic vectors of $S$-matchings. Since $|S|$ is odd, for an $S$-matching $M$ some edges of $M$ necessarily go outside $S$. Therefore the characteristic vector $\chi_M$ satisfies $\sum_{e\in C(S,S')} (\chi_M)_e \ge 1$. Hence any vector $x\in P_S$ satisfies $\sum_{e\in C(S,S')} x_e\ge 1$. In particular, for $x\in tP_S$, $\sum_{e\in C(S,S')} x_e\ge t$.\end{proof}
\noindent\textbf{Acknowledgements}

The author would like to thank Professor Martin Henk for helpful discussion and encouragement while the author was staying at Magdeburg university and studying on the topics in lattice points enumeration.

\end{document}